\documentclass[11pt]{article}
\usepackage{lscape}
 \usepackage[ansinew]{inputenc}
 \usepackage[dvips]{graphicx}
\usepackage[psamsfonts]{amssymb} \usepackage[all]{xy}  \usepackage[all]{xy}
 \usepackage[psamsfonts]{eucal}
 \usepackage[psamsfonts]{eucal}
\usepackage{amssymb, amsmath}
\usepackage{amsthm}
\usepackage{geometry}
\usepackage{enumerate}
\usepackage{float}
\newcommand{\be}{\begin{equation}}
\newcommand{\ee}{\end{equation}}

\newcommand{\bd}{\begin{displaymath}}
\newcommand{\ed}{\end{displaymath}}
\newcommand{\ben}{\begin{enumerate}}
\newcommand{\een}{\end{enumerate}}
\newtheorem{Theorem}{Theorem}

\title{Aspherical completions and rationally inert elements}
\author{Y. F\'elix and S. Halperin}
\begin{document}
\maketitle

\begin{abstract} Let $X$ be a connected space. An element $[f]\in \pi_n(X)$ is called rationally inert if
 $\pi_*(X)\otimes \mathbb Q \to \pi_*(X\cup_fD^{n+1})\otimes \mathbb Q$ is surjective. We extend the results of \cite{HL}  and prove in particular that if $X\cup_fD^{n+1}$ is a Poincar\'e duality complex and the algebra $H(X)$ requires at least two generators  then $[f]\in \pi_n(X)$ is rationally inert. On the other hand, if $X$ is rationally a wedge of at least two spheres and $f$ is rationally non trivial, then $f$ is rationally inert. Finally if $f$ is rationally inert then the rational homotopy of the homotopy fibre of the injection $X \to X\cup_fD^{n+1}$ is the completion of a free Lie algebra.
 \end{abstract}

\vspace{5mm}\noindent 
{\bf 2010MSC:} 55P62, 55P05

\vspace{1cm}
 In \cite{An2} and \cite{HL} the authors define and establish the properties of rationally inert elements in the homotopy groups of simply connected CW complexes $X$ of finite type: $[f]\in \pi_n(X)$ is \emph{rationally inert} if
 $$\pi_*(X)\otimes \mathbb Q \to \pi_*(X\cup_fD^{n+1})\otimes \mathbb Q$$
 is surjective. Our objective here is to use Sullivan completions $X\to X_{\mathbb Q}$ to extend the definitions to $[f]\in \pi_n(X)$, $n\geq 1$, where $X$ is any connected CW complex, and then to extend the principal results of \cite{HL} to this more general setting and establish several applications. For details about Sullivan completions the reader is referred to \cite{completions}
 
Inverse homotopy equivalences between the homotopy categories of connected CW complexes, $X$, and connected simplicial sets, $S$, are provided by $X\mapsto$ Sing$\,X$, the singular simplices in $X$, and by $S\mapsto \vert S\vert$, its Milnor realization. These identify a map $X\to \vert S\vert$ with a morphism Sing$\, X\to S$. For simplicity we denote both by
 $$X\to S,$$
 and refer to either a connected CW complex or a connected simplicial set simply as a \emph{connected space}.
 
 Additionally, for simplicity, we adopt the 
 
 \vspace{3mm} \noindent {\bf Convention.} Our base field is $\mathbb Q$. When the meaning is clear, we will suppress the differentials from the notation. For simplicity, we will also write
 $$ (-)^\vee := \mbox{Hom}(-, \mathbb Q), \hspace{5mm}\mbox{and } H(-):= H^*(-;\mathbb Q),$$
 for singular cohomology. Moreover, where there is no ambiguity we suppress the differential from the notation for a complex, and write $A$ instead of $(A,d)$. 
 
 \vspace{3mm}
 As detailed in \S 1 below, a Sullivan completion $X_{\mathbb Q}$ appears naturally as a simplicial set. Sullivan models and Sullivan completions are reviewed in \S 1. In particular,   if $X$ is simply connected and of finite type, then \cite[Theorem 15.11]{FHTI} its Sullivan completion induces an isomorphism $\pi_*(X)\otimes \mathbb Q \stackrel{\cong}{\longrightarrow} \pi_*(X_{\mathbb Q})$. Thus we extend the definition of rationally inert elements as follows:
 
 \vspace{3mm}\noindent {\bf Definition.} If $X$ is a connected space then $[f]\in \pi_n(X)$, some $n\geq 1$, is \emph{rationally inert} if the inclusion $i : X\to X\cup_fD^{n+1}$
 induces a surjection,
 $$\pi_*(i_{\mathbb Q}): \pi_*(X_{\mathbb Q}) \to \pi_*((X\cup_f D^{n+1})_{\mathbb Q}).$$
 
 \vspace{3mm} This condition can be characterized in   terms of the homotopy type of the fibre $F(f)$ of $i_{\mathbb Q}$ (Theorem 1). Applications are then provided in Theorems 2, 3 and 4. To state Theorem 1 we need the

 \vspace{3mm}\noindent {\bf Definition.} A connected space $Y$ is \emph{rationally wedge-like} if for some non-void linearly ordered set $S = \{\sigma\}$, and integers $n_\sigma >0$, there is a homotopy equivalence,
 $$Y\stackrel{\simeq}{\longrightarrow} \varprojlim_{\sigma_1<\dots <\sigma_r} \, (S^{n_{\sigma_1}}\vee \dots \vee S^{n_{\sigma_r}})_{\mathbb Q},$$
 where the inverse system is defined by the projections of $S^{n_{\sigma_1}}\vee \dots \vee S^{n_{\sigma_r}}$ on the sub wedges.
 
\vspace{3mm}\noindent {\bf Remark:} Note that in general $(X\vee Y)_{\mathbb Q}$ is different from $X_{\mathbb Q}\vee Y_{\mathbb Q}$ !

 \begin{Theorem} For any connected space $X$, a homotopy class, $[f]\in \pi_n(X)$, some $n\geq 1$, is rationally inert if and only if   the homotopy fibre $F(f)$ of $X_{\mathbb Q}\to (X\cup_fD^{n+1})_{\mathbb Q}$ is    rationally wedge-like.\end{Theorem}
\noindent Applications are then provided in Theorems 2, 3, 4 and 5.

\vspace{3mm} Theorem 3.3 in \cite{HL} is   a special case of Theorem 1 since in that case the homotopy fibre of $X\to X\cup_fD^{n+1}$ is rationally a wedge of spheres if and only if its rationalization is rationally wedge-like.

An example of rationally inert elements is provided by the following theorem, established for simply connected spaces in (\cite[Theorem 5.1]{HL}).

\begin{Theorem} If $X\cup_fD^{n+1}$ is a Poincar\'e duality complex and the algebra $H(X)$ requires at least two generators  then $[f]\in \pi_n(X)$ is rationally inert.
\end{Theorem}

\vspace{3mm}
 As described above, and in detail in (\cite[\S 1]{completions}) the Sullivan completion $X_{\mathbb Q}$ of a space $X$ is a simplicial set $\langle \land W\rangle$ constructed from a minimal Sullivan model for $X$. This is used in (\cite[\S 4]{completions}) to construct a completion, $\overline{H(\Omega X)}$, of the rational loop space homology of $X$. The homotopy fibre $F(f)$ of Theorem 1 also has the form $\langle \land Z\rangle$ for some minimal Sullivan algebra $\land Z$ (\S 3), although $\land Z$ may not be the Sullivan model of a space. Nevertheless, (\cite[\S 4]{completions}), for any minimal Sullivan algebra, $\land Z$, $\pi_*\Omega \langle \land Z\rangle$ is naturally a graded Lie algebra, complete with respect to a natural filtration. Its Lie bracket is given explicitly in terms of the Whitehead products in $\pi_*\langle \land Z\rangle$. We generalize (\cite[Theorem 3.3 (I)]{HL}) in 
 
\begin{Theorem} Suppose $X$ is a connected space and $[f]\in \pi_n(X)$, some $n\geq 1$, is rationally inert. Then $\pi_*(\Omega F(f))$ is the completion of a free sub Lie algebra, freely generated by a subspace $S\cong \overline{H_*(\Omega (X\cup_fD^{n+1})}$. \end{Theorem}

 \vspace{3mm} A general question   asks  what conditions on a group $G$ imply that    $ (BG)_{\mathbb Q}$ is aspherical; i.e., a $K(\pi, 1)$. This is true when $G$ is a finitely generated free group, when $G$ is the fundamental group of a Riemann surface or when $G$ is a right-angled Artin group (\cite{PS}, \cite{FH}). We consider here the one-relator groups, $\pi_1(X\cup_fD^2)$,   obtained  
 by adding a 2-cell to a   wedge of circles along a continuous map   $f: S^1\to X$.  The well known Lyndon theorem (\cite{Ly}, \cite{Pu},\cite{DV}) states that if   $f$ is not a proper power, then $X\cup_fD^2$ is aspherical.  In general it may happen that a connected space $X$ is aspherical, but $X_{\mathbb Q}$ is not. However, the spaces considered by Lyndon remain aspherical when rationalized:

\begin{Theorem} If $X$ is a wedge of at least two circles then any non zero $[f]\in \pi_1(X)$ is rationally inert;   equivalently, $(X\cup_fD^2)_{\mathbb Q}$ is aspherical. \end{Theorem}

 \vspace{3mm}\noindent {\bf Remark.} Note that even if $f$ is a proper power, where Lyndon's theorem does not apply, it is true that $(X\cup_fD^2)_{\mathbb Q}$ is aspherical.

 \vspace{3mm} Finally recall a famous unsolved problem of JHC Whitehead \cite{W}: is a subcomplex of an aspherical two-dimensional CW complex aspherical ? As observed by Anick \cite{An} it is sufficient to consider the case that both subcomplexes share the same $1$-skeleton and base point. The problem then reduces to the question: If $X$ is a finite 2-dimensional connected CW complex and $X\cup \left( \amalg_{k=1}^p D_k^2\right)$ is aspherical, is $X$ aspherical ?
 
 In \cite{An} Anick provides a positive answer to an analogous question for simply connected rational spaces. Here we have a positive answer for   Sullivan completions of   connected spaces.
 
\begin{Theorem} If $X$ is a     connected space and $\left( X\cup \amalg_{k=1}^p D^2_k\right)_{\mathbb Q}$ is aspherical, then $X_\mathbb Q$ is aspherical.\end{Theorem}

 \section{Sullivan models and Sullivan completions}

\vspace{3mm}
 We review briefly the basic facts and notation from Sullivan's theory. For details the reader is referred to   \cite{completions}. A \emph{$\Lambda$-algebra} is a commutative differential graded algebra (cdga) of the form $(\land V,d)$, where $V= V^{\geq 0}$ is a graded vector space and $\land V$ is the free graded commutative algebra generated by $V$. Moreover the differential is required to satisfy the \emph{Sullivan condition}: $V = \cup_{n\geq 0} V(n)$, where
 $$V(0) = V\cap \mbox{ker}\, d \hspace{5mm}\mbox{and } V(n+1) = V\cap d^{-1}(\land V(n)).$$
 Here $V$ is a \emph{generating vector space} for $\land V$. If $V= V^{\geq 1}$ then $V$ is a \emph{Sullivan algebra}.

Moreover, $\land V = \oplus_{p\geq 0} \land^pV$, where $\land^pV$ denotes the linear span of the monomials in $V$ of length $p$; $p$ is called the \emph{wedge degree}. In particular, a $\Lambda$-algebra is \emph{minimal} if $d : V\to \land^{\geq 2} V$ and \emph{quadratic} if $d : V\to \land^2V$. Thus associated with a minimal $\Lambda$-algebra $(\land V,d)$ is the quadratic $\Lambda$-algebra $(\land V,d_1)$ defined by: $d_1v$ is the component of $dv$ in $\land^2V$.
 
  Note that if $V= V^{\geq 1}$, then the inclusion of a  subspace $W\subset \land^{\geq 1} V$ extends to an isomorphism $\land W\stackrel{\cong}{\to} \land V$ if and only if $W\oplus \land^{\geq 2}V = \land V$. In this case $\land W$ satisfies the same condition as $\land V$: the definition of a Sullivan algebra does not depend on the choice of generating vector space. Observe as well that if $V= V^{\geq 1}$ then the natural map
 \begin{eqnarray}
 \xymatrix{\land V\ar[r]^\cong & \prod_p\land^pV}
 \end{eqnarray}
 is an isomorphism.
 
With each connected space $Y$ is associated a cdga $A_{PL}(Y)$ and a unique isomorphism class of minimal Sullivan algebras $(\land V,d)$ characterized by the existence of a quasi-isomorphism $(\land V,d)\stackrel{\simeq}{\to} A_{PL}(Y)$. By definition $(\land V,d)$ is the \emph{minimal Sullivan model} of $Y$. Among their properties are the natural isomorphisms $H(\land V,d)\cong H(Y)$ of graded algebras. Moreover, any map, $f: X\to Y$   determines a "homotopy class" of morphisms, $\varphi: \land V\to \land W$, from the minimal Sullivan model of $Y$ to that of $X$; $\varphi$ is a \emph{Sullivan representative} of $f$.

  On the other hand, the construction of Sullivan completions is accomplished   by a functor   associating to a $\Lambda$-algebra, $\land W$,  a simplicial set $\langle \land W\rangle$, with the property that \emph{$<\, \,>$ converts direct limits to inverse limits}. In particular, if 
  $\land W$ is a minimal Sullivan model of a connected space $X$ then this determines a based homotopy class of maps  
 $$  X \longrightarrow X_{\mathbb Q}:= \langle \land W\rangle,$$
  the \emph{Sullivan completion} of $X$. In particular, 
 if $\varphi : \land V\to \land W$ is a Sullivan representative of $f: X\to Y$ then
$$f_{\mathbb Q}= \langle \varphi\rangle : \langle \land W\rangle \to \langle \land V\rangle.$$

Moreover, (\cite[Theorem 1.3]{RHTII}) for any minimal Sullivan algebra, $\land W$, there is a natural bijection $\pi_*(\langle \land W\rangle)\cong \mbox{Hom}(\land^{\geq 1}W/\land^{\geq 2}W)$, and the isomorphism $W\stackrel{\cong}{\to} \land^{\geq 1}W/\land^{\geq 2}W$ then induces a bijection
$$\pi_*(\langle \land W\rangle) \cong  W^\vee.$$
Therefore, for any morphism $\varphi : \land V\to \land W$ of minimal Sullivan algebras, it follows that $\pi_*(\langle \varphi\rangle)$ is surjective if and only if $\varphi : \land^{\geq 1}V/\land^{\geq 2}V \to \land^{\geq 1}W/\land^{\geq 2}W$ is injective, or equivalently, if   the generating vector space $W\subset \land W$ can be chosen so that $\varphi : V\to W$ is the inclusion of a subspace. In this case
$$\pi_*(\langle \varphi\rangle) =  \varphi^\vee :  W^\vee\to  V^\vee.$$

Now a general morphism $\varphi : \land V\to \land W$ of Sullivan algebras factors (\cite[Theorem 3.1]{RHTII}) as
$$\xymatrix{\land V \ar[r]^-\eta & \land V\otimes \land Z \ar[r]^-\gamma_-\simeq & \land W,}$$
in which (i) $\eta (v) = v\otimes 1$, (ii) $\gamma$ is a quasi-isomorphism, (iii) $Z= Z^{\geq 0}$,  (iv) $Z= \cup_n Z(n)$  satisfying
$$ Z(0)= Z\cap d^{-1}(\land V) \hspace{5mm}\mbox{and } Z(n+1) = Z\cap d^{-1} (\land V\otimes \land Z(n)),$$
and (v)   the quotient $(\land Z, \overline{d}) = \mathbb Q \otimes_{\land V}(\land V\otimes \land Z)$ is a minimal $\Lambda$-algebra. Here  
$\land V\otimes \land Z$ is a \emph{minimal $\Lambda$-extension} of $\land V$. 

\vspace{3mm}\noindent {\bf Remark.} If $\pi_*\langle \varphi\rangle$ is surjective we take $\eta= \varphi$ to be  an inclusion $V\to W$ and $\land W = \land V\otimes \land Z$.

 \vspace{2mm}
 In particular, with each minimal Sullivan algebra $(\land V,d)$ is associated a unique isomorphism class of $\Lambda$-extensions, $(\land V\otimes \land U,d)$, its \emph{acyclic closures}. These are characterized by the following two properties: (i) the augmentation $\land V\to \mathbb Q$ extends to a quasi-isomorphism $\land V\otimes \land U\stackrel{\simeq}{\to} \mathbb Q$ with $U\to 0$, and (ii) the quotient differential in $\land U= \mathbb Q\otimes_{\land V} (\land V\otimes \land U)$ is zero.

Finally, a minimal Sullivan algebra $\land V$ determines the graded \emph{homotopy Lie algebra} $L_V = (L_V)_{\geq 0}$ given by
 $$s(L_V)_p = \mbox{Hom}(V^{p+1})$$
 and
 $$<v, s[x,y]> = (-1)^{1+deg\, y} <d_1v, sx, sy>\,, \hspace{5mm} v\in V, x,y\in L_V.$$
(Here $s$ is the degree 1 \emph{suspension isomorphism}.) Thus
\begin{eqnarray}
\label{alpha}
s(L_V) = \pi_*\langle \land V\rangle.
\end{eqnarray}

 \section{Rationally wedge-like spaces}

\vspace{3mm}\noindent {\bf Lemma 1.}  {\sl The following two conditions on a minimal Sullivan algebra, $\land Z$, are equivalent:
\begin{enumerate}
\item[(i)] The generating vector space $Z\subset \land Z$ can be chosen so that 
$$Z\cap \mbox{ker}\, d\stackrel{\cong}{\longrightarrow} H^{\geq 1}(\land Z).$$
\item[(ii)] $\land Z$ is the minimal Sullivan model of a cdga $A = \mathbb Q \oplus A^{\geq 1}$ in which the differential and products in $A^{\geq 1}$ are zero.
\end{enumerate}
If these hold then $Z$ can be chosen so that $Z\cap \mbox{ker}\, d \stackrel{\cong}{\to} H^{\geq 1}(\land Z)$ and $(\land Z,d)$ is quadratic.}
 
\vspace{3mm}\noindent {\sl proof:} If (i) holds let $A$ be the quotient of $\land Z$ by $\land^{\geq 2}Z$ and by a direct summand of the image of ker$\,d$ in $Z$. If (ii) holds set $V_0 = A^{\geq 1}$ and define a quadratic Sullivan algebra $\land V$ by setting $V(k) = \oplus_{j\leq k} V_k$, with $d : V_{k+1}\stackrel{}{\to} \land^2V(k) \cap \mbox{ker}\, d$ inducing an isomorphism in homology. Then $(\land V,d)$ has zero homology in wedge degree $2$, and it follows that $\land V$ has zero  homology in wedge degrees $\geq 2$. Hence $\land V$ is a quadratic Sullivan model for $A$. Thus  $\land V \cong \land Z$, and so $Z$ can be chosen so that $d : Z\to \land^2Z$. Thus the final assertion is part of (\cite[Proposition 6]{completions}). \hfill$\square$

\vspace{3mm}\noindent {\bf Example:} Finite wedges of spheres: $S = S^{\sigma_1}\vee \dots \vee S^{\sigma_k}$.

\vspace{2mm}The quasi-isomorphism $A_{PL}(S)\to \oplus_{\mathbb Q} A_{PL}(S^{\sigma_i}) \simeq \oplus_{\mathbb Q} H(S^{\sigma_i})$ 
identifies the minimal Sullivan model of $S$ as a minimal Sullivan algebra $\land Z$ satisfying the conditions of Lemma 1. Here $Z\cap \mbox{ker}\, d$ has a basis $z_1, \dots , z_k$ representing orientation classes of $S^{\sigma_1}, \dots , S^{\sigma_k}$.

Now choose elements $x_i $ in the homotopy Lie algebra $L_S$ of $S$ so that $\langle z_i, sx_j\rangle = \delta_{ij}$.  
  The $x_j$ then freely generate a free sub Lie algebra $E\subset L_S$. In fact,   the rescaling argument in (\cite[p.230]{RHTII}) generalizes to reduce to the case $S= S^{\geq 2}$, in which case the result is established in \cite[\S 23, Example 2]{FHTI}. Moreover, 
 it follows from \cite[Chap. 2]{RHTII} that 
 $$L_S = \varprojlim_n L_S/L_S^n$$
 where $L_S^n$ is the ideal spanned by the iterated commutators in $L_S$ of length $n$. According to \cite[Chapter 2]{RHTII}, the $x_{\sigma_i}$ map to a basis of $L_S/L_S^2$ and hence the inclusion $E\hookrightarrow L_S$ induces isomorphisms $E/E^n \stackrel{\cong}{\to} L_S/L_S^n$.

\vspace{3mm}\noindent {\bf Proposition 1.} {\sl A connected space $F$ is rationally wedge-like if and only il it has the form $F = \langle \land Z\rangle$, where $\land Z$ satisfies the equivalent conditions of Lemma 1.}

\vspace{3mm}\noindent {\sl proof:}  Suppose first that $F = \langle \land Z\rangle$, where $\land Z$ satisfies the conditions of Lemma 1, and pick  a linearly ordered basis of $Z\cap \mbox{ker}\, d$. 
Then each finite subset $z_{\sigma_1}<\dots <z_{\sigma_k}$ determines an inclusion
 $$\land Z(\sigma_1, \dots , \sigma_k) \hookrightarrow \land Z$$
 of quadratic Sullivan algebras  with $Z(\sigma_1, \dots , \sigma_k)\subset Z$, and for which $\{z_{\sigma_i}\}$ is a basis of $H^{\geq 1}(\land Z(\sigma_1, \dots , \sigma_p))$, and    $\land Z(\sigma_1, \dots , \sigma_k)$ is   a Sullivan model for $S^{\sigma_1}\vee \dots \vee S^{\sigma_k}$. Moreover,    the inclusions $\land Z(\sigma_{i_1}, \dots \sigma_{i_r}) \to \land Z(\sigma_1, \dots , \sigma_k)$ are Sullivan representatives for the projections $S^{\sigma_1}\vee \dots \vee S^{\sigma_k}\to S^{\sigma_{i_1}} \vee \dots \vee S^{\sigma_{i_r}}.$
 
  Now $$\land Z = \varinjlim_{\sigma_1<\dots < \sigma_k} \, \land Z(\sigma_1, \dots , \sigma_k)$$
 and so
 $$\langle \land Z\rangle = \varprojlim_{\sigma_1<\dots < \sigma_k} \langle \land Z (\sigma_1, \dots , \sigma_k)\rangle = \varprojlim_{\sigma_1<\dots <\sigma_k} \, (S^{\sigma_1} \vee \dots \vee S^{\sigma_k})_{\mathbb Q}.$$

 In the reverse direction, suppose $F$ is rationally wedge like, so that
 $$F = \varprojlim_{\sigma_1<\dots < \sigma_k} (S^{n_{\sigma_1}} \vee \dots \vee S^{n_{\sigma_k}})_{\mathbb Q}.$$
 Then let $\land Z$ be a Sullivan algebra satisfying the conditions of Lemma 1 in which $Z \cap \mbox{ker}\, d$ has a basis $\{z_{\sigma}\}$ of degrees $n_{\sigma}$. Thus any subset $\sigma_{1} <\dots <\sigma_{k}$ determines a sub Sullivan algebra $Z(\sigma_{1}\dots \sigma_{k})$ by the requirement that $Z(\sigma_{1}\dots \sigma_{k})= \cup_n\, Z(\sigma_{1}\dots \sigma_{k};n)$ in which
 $$Z(\sigma_{1}\dots \sigma_{k};0) = \oplus_i \mathbb Q z_{\sigma_i}$$
 and
 $$Z(\sigma_{1}\dots \sigma_{k};n+1) = Z \cap d^{-1} (\land^2Z(\sigma_{1}\dots \sigma_{k};(n))).$$
 This gives as above that
 $$\langle \land Z\rangle = \varprojlim_{\sigma_{1}<\dots <\sigma_{k}}\,\, \langle \land Z(\sigma_{i}\dots \sigma_{k})\rangle = \varprojlim_{\sigma_{1}<\dots <\sigma_{k}}\, \, \left( S^{n_{\sigma_{1}}} \vee \dots \vee S^{n_{\sigma_{k}}} \right)_{\mathbb Q} = F.$$
 
 \hfill$\square$

\vspace{3mm}\noindent {\bf Corollary 1.}  \emph{If $X = \vee_\sigma S^{n_\sigma}$ is a wedge of spheres, then $X_{\mathbb Q}$ is rationally wedge-like. If all the spheres are circles then $X_{\mathbb Q}$ is aspherical.}

\vspace{3mm}\noindent {\bf Corollary 2.} \emph{If $\langle \land Z\rangle$ is rationally wedge-like and dim$\, H^{\geq 1}(\land Z)>1$, then the sum of the solvable ideals in $L_Z$ is zero.}

\vspace{3mm}\noindent {\sl proof:} It follows from Lemma 1 that cat$(\land Z)= 1$, and so from \cite{depthLS}, Sdepth$\, L_Z<\infty$. Now \cite[Theorem 1]{Sdepth} asserts that the sum, rad$\, L_Z$, of the solvable ideals in $L_Z$ is finite dimensional, and that $L_Z$ acts nilpotently in rad$\, L_Z$. In particular, if rad$\, L_Z\neq 0$ then the center of $L_Z$ is non-zero. Let $x\in L_Z$ be an element in the center.

Since $\langle \land Z\rangle$ is rationally wedge-like, $\land Z= \varinjlim_{\mathcal S} \land Z(\sigma_1, \dots , \sigma_k)$ where $\land Z(\sigma_1, \dots , \sigma_k)$ is the minimal Sullivan model of a wedge of $k$ spheres, and ${\mathcal S}$ has by hypothesis at least two elements. Then $L_Z= \varprojlim_{\mathcal S} L_{Z(\sigma_1, \dots , \sigma_k)}$, and the maps $L_Z\to L_{Z(\sigma_1, \dots , \sigma_k)}$ are surjective. Thus if $x\neq 0$ it maps to a non-zero element in some $L_{Z(\sigma_1, \dots , \sigma_k)}$ with $k> 1$. This would contradict the Example above. \hfill$\square$

\vspace{3mm}\noindent {\bf Remark.} Rationally wedge-like spaces provide examples of minimal Sullivan algebras $\land Z$ for which $\langle \land Z\rangle$ is not the Sullivan completion of a space. For example, suppose $Z= Z^3$ has a countably infinite basis, so that $\pi_*\langle \land Z\rangle = \pi_3\langle \land Z\rangle= (Z^3)^\vee$.

Thus for any minimal Sullivan algebra $\land V$, the condition $\langle \land V\rangle = \langle \land Z\rangle$ would imply that $V= V^3$ and $(V^3)^\vee \cong (Z^3)^\vee$. But if $\land V$ were the minimal model of a space $X$ then we would have $V^3\cong H^3(X)= H_3(X)^\vee$ and so either dim$\, V^3<\infty$ or card$\, (V^3)\geq $ card$\,  \mathbb R$. In the second case, card$\, ((V^3)^\vee)> $ card$\, \mathbb R$ and so $(V^3)^\vee$ and $(Z^3)^\vee$ are not isomorphic.

\vspace{3mm}\noindent {\bf Proposition 2.} {\sl Suppose $X$ and $Y$ are connected spaces, one of which has rational homology of finite type. Then
\begin{enumerate}
\item[(i)] The homotopy fibre, $F$, of the natural map
$$i_{\mathbb Q}: (X\vee Y)_{\mathbb Q}\to (X\times Y)_{\mathbb Q}$$
is rationally wedge-like.
\item[(ii)] If $X_{\mathbb  Q}$ and $Y_{\mathbb Q}$ are aspherical then so are $F$ and $(X\vee Y)_{\mathbb Q}$.
\end{enumerate}}

\vspace{3mm} This result is analogous to the fact that the usual fibre of the injection $X\vee Y\to X\times Y$ is the join of $\Omega X$ and $\Omega Y$ and thus a suspension. (But note that $(X\vee Y)_{\mathbb Q}$ may be different from $X_{\mathbb Q}\vee Y_{\mathbb Q}$.)

Proposition 2   follows easily from a result about Sullivan algebras (Proposition 3, below). For this, consider minimal Sullivan algebras,   $\land W$ and $\land Q$. The natural surjection $\land W\otimes \land Q \to \land W\times_{\mathbb Q}\land Q$ is surjective in homology, and so extends to a minimal Sullivan model
$$\varphi : \land T:= \land W\otimes \land Q\otimes \land R \stackrel{\simeq}{\longrightarrow} \land W\times_{\mathbb Q}\land Q.$$
Filtering by wedge degree then yields a morphism
$$\varphi_1 : (\land T, d_1) \to (\land W, d_1)\times_{\mathbb Q}(\land Q, d_1)$$
between the associated bigraded cdga's. (Here $(\land -,d_1)$ is the associated quadratic Sullivan algebra.) 

\vspace{3mm}\noindent {\bf Proposition 3.} {\sl With the hypotheses and notation above,
\begin{enumerate}
\item[(i)] $\langle \land R\rangle$ is rationally wedge-like.
\item[(ii)] $\varphi_1$ is a quasi-isomorphism.
\end{enumerate}}

\vspace{3mm}\noindent {\sl proof:} (i) Let $\land W\otimes \land U_W$ and $\land Q\otimes \land U_Q$ denote the respective acyclic closures. 
Then  $\land R$ is quasi-isomorphic to
 $$\land T \otimes_{\land W\otimes \land Q} \land W\otimes \land U_W\otimes \land Q\otimes \land U_Q \simeq A:= (\land W\oplus_{\mathbb Q}\land Q) \otimes \land U_W\otimes \land U_Q.$$
Dividing $A$ by the ideal generated by $W$ yields the short exact sequence
 $$0\to \land^{\geq 1}W\otimes \land U_W\otimes \land U_Q \to A \to  \land Q\otimes \land U_W\otimes \land U_Q \to 0.$$
 
 Decompose the differential in $\land W\otimes \land U_W$ in the form $d= d_1+d'$ with $d_1(W)\subset \land^2W$, $d_1(U_W)\subset W\otimes \land U_W$, $d'(W)\subset \land^{\geq 3}W$ and $d'(U_W)\subset \land^{\geq 2}W\otimes \land U_W$. Then $d_1$ is a differential and $(\land W\otimes \land U_W,d_1)$ is the acyclic closure of $(\land W,d_1)$. Choose a direct summand,  $S$, of $d_1(\land^{\geq 1}U_W) $ in $W\otimes \land U_W$. Then $I = (\land^{\geq 2}W\otimes \land U_W)\oplus S$ is acyclic for the differential $d_1$ and therefore also for the differential $d$. 
Thus $J = I\otimes \land U_Q$ is an acyclic ideal in $A$ and $A\stackrel{\cong}{\to} A/J$.  

Now consider the short exact sequence
$$0\to (\land ^{\geq 1}W\otimes \land U_W\otimes \land U_Q)/J \to A/J\to \land Q\otimes \land U_W\otimes \land U_Q\to 0.$$
  The inclusion of $\land U_W$ in the right hand term is a quasi-isomorphism. This yields a quasi-isomorphism
 $$  d_1(\land^{\geq 1}U_W) \otimes \land^{\geq 1}U_Q\simeq A/J\simeq A.$$
 Since  the   differential and the multiplication in $d_1(\land^{\geq 1}U_W) \otimes \land^{\geq 1}U_Q$ are zero, it follows from Proposition 1   that $\langle \land R\rangle$ is rationally wedge-like.

\vspace{3mm} (ii) The surjection $(\land W\otimes \land Q,d_1)\to (\land W\times_{\mathbb Q}\land Q,d_1)$ extends to a quasi-isomorphism
$$\widehat{\varphi} : \land \widehat{T} := (
\land W\otimes \land Q\otimes \land \widehat{R}, \delta) \to (\land W\times_{\mathbb Q} \land Q,d_1)$$
from a minimal Sullivan algebra. We first show that $\widehat{R}$ can be chosen so that $(\land \widehat{T}, \delta)$ is quadratic. Then we extend $\delta$ to a differential $\widehat{d}= \sum_{i\geq 1} \widehat{d}_i$ in which $\widehat{d}_1=\delta$ and 
$$\widehat{d}_i : \widehat{T}\to \land^{i+1}\widehat{T} \hspace{3mm}\mbox{and } \widehat{\varphi}\circ \widehat{d} = d\circ \varphi.$$
It is automatic that $(\land \widehat{T}, \widehat{d})$ will be a minimal Sullivan algebra. Moreover, filtering by wedge degree shows that $\widehat{\varphi}$ is a quasi-isomorphism and so $\land \widehat{T}$ is a minimal Sullivan model for $\land W\times_{\mathbb Q}\land Q$. In particular this identifies $\widehat{T}$ with $T$, $R$ with $\widehat{R}$  and $\widehat{\varphi}$ with $\varphi$, thereby establishing (ii).

\vspace{2mm}To accomplish the first step, define
 $d_1 : U_W\to W\otimes \land U_W$   and   $d_1 : U_Q\to Q\otimes \land U_Q$ as in (i). Assign  $\land W$ and $\land Q$ wedge degree as a second degree and assign $U_W$ and $U_Q$   second degree $0$. Then $(\land W\otimes \land U_W,d_1)$ and $(\land Q\otimes \land U_Q,d_1)$ are the respective acyclic closures of $(\land W,d_1)$ and $(\land Q, d_1)$, and $d_1$ increases the second degree by 1. Now $\widehat{\varphi}$ and $\widehat{T}$ 
 may be constructed so that $ \widehat{R}$ is equipped with a second gradation for which $\delta$ increases the second degree by one and $\widehat{\varphi}$ is bihomogeneous of degree zero. 

The argument in the proof of (i) now yields a sequence of bihomogeneous quasi-isomorphisms connecting
$$\mathbb Q \oplus \left( d_1(\land^+U_W)\otimes \land^+U_Q\right)\simeq \land \widehat{R}.$$
Thus $H^{\geq 1}(\land \widehat{R})$ is concentrated in second degree 1. Therefore $\land \widehat{R}$ satisfies condition (i) of Proposition 1, and it follows that we may choose $\widehat{R}$ so that the quotient cdga $\land \widehat{R}$ is quadratic and $H^{\geq 1}(\land \widehat{R})$ embeds in $\widehat{R}$. 
This implies that $\widehat{R}$ is concentrated in second degree 1 and that 
$$\delta: \widehat{R}\to     \land^2(W\oplus Q\oplus T).$$
In particular, $(\land \widehat{T},\delta)$ is a quadratic Sullivan algebra.

\vspace{3mm} The construction of $\widehat{d}$ proceeds as follows. Write the differential in $\land W\times_{\mathbb Q}\land Q$ as $d= \sum_{r\geq 1} d_r$ in which $d_r$ is a derivation raising wedge degree by $r+1$. Thus for each $r$, $\sum_{i+j=r} d_id_j = 0$. Now we construct by induction a sequence of derivations $\widehat{d}_1=\delta, \dots , \widehat{d}_r\dots$, in $\land \widehat{T}$ in which $\widehat{d}_r$ increases the wedge degree by $r+1$, and 
$$\sum_{i+j=r}\widehat{d}_i\widehat{d}_j = 0\hspace{3mm}\mbox{and } \widehat{\varphi}\widehat{d}_i = d_i\widehat{\varphi}.$$
Thus, in view of (i),  $\widehat{d}:= \sum \widehat{d}_i$ will define a differential in $\widehat{T}$,  $(\land \widehat{T}, \widehat{d})$ will be a Sullivan algebra, and $$\varphi : (\land \widehat{T}, \widehat{d}) \to (\land W\times_{\mathbb Q}\land Q,d)$$ will be a cdga morphism. Filtering by wedge degree shows that $\widehat{\varphi}$ is a quasi-isomorphism. 

It remains to construct the $\widehat{d}_i$, $i\geq 2$. For this, 
 set $\widehat{T}(k) = W\oplus Q \oplus \widehat{R}^{\leq k}$. Since $(\land \widehat{T}, \delta)$ is a Sullivan algebra it follows that each $\widehat{R}^k$ is the union of an increasing family of subspaces $F^p(\widehat{R}^k)$ such that 
$$\delta : F^0(\widehat{R}^k) \to \land \widehat{T}(k-1) \hspace{3mm}\mbox{and } \delta : F^{p+1}(\widehat{R}^k) \to \land \widehat{T}(k-1) \otimes \land F^p(\widehat{R}^k).$$
Set $\widehat{d}_1=\delta$ and assume by induction that $\widehat{d}_1, \dots , \widehat{d}_r$ have been constructed, and that $\widehat{d}_{r+1}$ has been constructed in $\widehat{R}^{<k}\oplus F^p(\widehat{R}^k)$.

Let $y_i$ be a basis for a direct summand of $F^p(\widehat{R}^k)$ in $F^{p+1}(\widehat{R}^{k})$. Then
$$\renewcommand{\arraystretch}{1.3}\begin{array}{ll}
\widehat{\varphi} (\widehat{d}_{r+1}\widehat{d}_1y_i) &= \widehat{d}_{r+1}\widehat{d}_1\psi y_i = -\widehat{d}_1\widehat{d}_{r+1}\widehat{\varphi}y_i - \displaystyle\sum_{j=2}^r (\widehat{d}_j\widehat{d}_{r+2-j})\widehat{\varphi} y_i
\\ &= - \widehat{d}_1\widehat{d}_{r+1}\widehat{\varphi} y_i - \displaystyle\sum_{j=2}^r \widehat{\varphi} (\widehat{d}_j\widehat{d}_{i+2-j})y_i.
\end{array}
\renewcommand{\arraystretch}{1}
$$

It follows that 
$$\widehat{d}_1\widehat{\varphi} \, (\, \widehat{d}_{r+1}\widehat{d}_1 + \sum_{j=2}^r \widehat{d}_j\widehat{d}_{r+2-j}\,)\,y_i = 0.$$
Since $\widehat{\varphi}$ is a surjective quasi-isomorphism with respect to $\widehat{d}_1$ and $d_1$, this implies that 
$$(\widehat{d}_{r+1}\widehat{d}_1 + \sum_{j=2}^r \widehat{d}_j\widehat{d}_{r+2-j})y_i = \widehat{d}_1\Phi_i$$
with $\widehat{\varphi} \Phi_i = -\widehat{d}_{r+1}\widehat{\varphi} y_i$. Extend $\widehat{d}_{r+1}$ to $F^{p+1}(\widehat{R}^k)$ by setting $\widehat{d}_{r+1}y_i = - \Phi_i$. 

\hfill$\square$

\vspace{3mm}\noindent {\sl proof of Proposition 2:} (i) Let $\land W$ and $\land Q$ be the minimal Sullivan models of $X$ and $Y$. A Sullivan representative of the inclusion $i : X\vee Y\to X\times Y$ is then the inclusion
$$\land W\otimes \land Q \to \land T:=\land W\otimes \land Q\otimes \land R.$$
It follows that $i_{\mathbb Q}$ is the surjection
$$\langle \land W\otimes \land Q\rangle \longleftarrow \langle \land W\otimes \land Q\otimes \land R\rangle.$$
But this surjection is a fibration (\cite[Proposition 17.9]{FHTI}) with fibre $\langle \land R\rangle$, which is a rationally wedge-like by Proposition 3.

 (ii) When $X_{\mathbb Q}$ and $Y_{\mathbb Q}$ are aspherical, then $U_W$ and $U_Q$ are concentrated in degree $0$ and $W$ is concentrated in degree 1. This shows that $F$ is aspherical. Since one of $X, Y$ has rational homology of finite type, $(X\times Y)_{\mathbb Q}= X_{\mathbb Q}\times Y_{\mathbb Q}$ is aspherical. We deduce then from the homotopy sequence of the fibration $F\to (X\vee Y)_{\mathbb Q}\to (X\times Y)_{\mathbb Q}$ that $(X\vee Y)_{\mathbb Q}$ is also aspherical.
 
\hfill$\square$

 \section{Cell attachments and Theorem 1}
 
 Before undertaking the proof of Theorem 1 we   set up the basic framework that translates the topology of a cell attachment to Sullivan's theory, and establish two preliminary Propositions.

 \vspace{3mm} Suppose $f : S^n\to X$ is the map of Theorem 1, and denote by $(\land W,d)$ the Sullivan minimal model of $X$. A Sullivan representative of $f$ is a morphism from $\land W$ to the minimal model of $S^n$. Composing with the quasi-isomorphism from that model to $H(S^n)$ gives a morphism  $\psi : \land W\to H(S^n)$. Now define  a linear map of degree $-n$,
 $$\varepsilon : \land W\to \mathbb Q,$$
 by setting $\varepsilon (1) = 0$ and $\psi (\Phi) = \varepsilon (\Phi)\cdot [S^n]$, $\Phi \in \land^{\geq 1}W$, where $[S^n]$ denotes an orientation class in $S^n$. In particular, $\varepsilon \circ d= 0$ and $\varepsilon (\land^{\geq 2}W)= 0$.

Now define a cdga $(\land W\oplus \mathbb Q a,D)$   as follows: deg$\, a= n+1$, $a^2=a\cdot \land^+W= 0$, and
 $$Da=0\hspace{5mm}\mbox{and } D\Phi= d\Phi+\varepsilon (\Phi)a, \hspace{1cm} \Phi\in \land W.$$
 By \cite[(13)b and (13)d]{FHTI}, division by $a$  yields the commutative diagram,
 \begin{eqnarray}
 \label{first}
 \xymatrix{\mathbb Q a\,\, \ar@{^{(}->}[r] & (\land W\oplus \mathbb Q a, D) \ar@{->>}[r] & (\land W,d)\\
 & (\land V,d), \ar[u]^\tau_\simeq \ar[ru]^\lambda}
 \end{eqnarray} 
 in which $(\land V,d)$ is a minimal Sullivan model for $X\cup_fD^{n+1}$, and $\lambda$ is a Sullivan representative for the inclusion $i : X\to X\cup_fD^{n+1}$.
 In particular, $i_{\mathbb Q}: X_{\mathbb Q}\to (X\cup_fD^{n+1})_{\mathbb Q}$
  is identified with $\langle \lambda\rangle : \langle \land W\rangle \to \langle \land V\rangle$.

\vspace{3mm}  As described in $\S 1$, 
$\lambda$ factors as
 $$\xymatrix{ \lambda: (\land V,d) \ar[r]^\eta & (\land V\otimes \land Z,d) \ar[r]^\gamma_\simeq & (\land W,d),}$$
 in which $\land V\otimes \land Z$ is a $\Lambda$-extension of $\land V$,   $\gamma $ is a quasi-isomorphism, and the quotient   $$(\land Z, \overline{d}):= \mathbb Q\otimes_{\land V}(\land V\otimes \land Z,d)$$
 is a minimal $\Lambda$-algebra. Since $H^1(i)$ is injective,   it follows that $\lambda : V^1\to W^1$ is injective. Therefore $Z= Z^{\geq 1}$ and $\land Z$ is a minimal Sullivan algebra. 
 
Further, because $\gamma$ is a quasi-isomorphism of Sullivan algebras, $\langle \gamma\rangle$ is a 
 homotopy equivalence, which (up to homotopy) identifies $\langle \eta\rangle$ with $\langle \lambda\rangle$. But (\cite[Proposition 17.9]{FHTI}) $\langle \eta\rangle$ is the projection of a Serre fibration with fibre $\langle \land Z\rangle$. Thus $\langle \land Z\rangle$, the homotopy fibre of $\langle \lambda\rangle$, and the homotopy fibre $F(f)$ of $i_{\mathbb Q}$, all have the same homotopy type:
 \begin{eqnarray}
 \label{dites3}
 \langle \land Z\rangle \simeq F(f).
 \end{eqnarray}
 
 On the other hand, we have

 \vspace{3mm}\noindent {\bf Proposition 4.} {\sl With the hypotheses and notation of (3), let $\land V\otimes \land U$ be the acyclic closure of $\land V$. Then   there is a degree 1 isomorphism,
 $$H^{\geq 1}(\land Z, \overline{d}) \stackrel{\cong}{\longrightarrow} \mathbb Q a \otimes \land U,$$
 and $H^{\geq 1}(\land Z, \overline{d})\cdot H^{\geq 1}(\land Z, \overline{d}) = 0.$}

 \vspace{3mm}\noindent {\sl proof:}
 First observe that in diagram (3), $\tau\Phi = \lambda \Phi+ \alpha (\Phi)a$. Thus $\tau$ must coincide with $\lambda$ in $\land^{\geq 2}V$, and that also $D\circ \tau = D\circ \lambda$. Thus for $\Phi\in \land V$,
 $$d(\lambda \Phi)+ \varepsilon (\lambda \Phi)a = D(\lambda \Phi) = D(\tau \Phi) = \tau d\Phi = \lambda (d\Phi) = d(\lambda \Phi).$$
Hence
 \begin{eqnarray}
 \label{4}
 \varepsilon \circ \lambda = 0.
 \end{eqnarray}
 
  Now let $\land V\otimes \land U$ be the acyclic closure of $\land V$. Apply $-\otimes_{\land V}\land V\otimes \land U$ to  diagram (\ref{first}) to obtain a short exact sequence of complexes,
 \begin{eqnarray}
 \label{nn}
0\to  \mathbb Q a \otimes \land U \to (\land W\oplus \mathbb Q a)\otimes_{\land V}(\land V\otimes \land U) \to \land W\otimes \land U\to 0,
 \end{eqnarray}
 in which the differential in $\mathbb Q a \otimes \land U$ is zero and the homology of the central complex is $\mathbb Q \, 1$ in positive degrees. It follows that $H^0(\land W\otimes \land U) = \mathbb Q 1$ and that the connecting homomorphism is an isomorphism of degree 1. By (\ref{4}), $\varepsilon$ vanishes   on $\land V$, and hence $(\varepsilon\otimes id)\circ (\lambda \otimes id)= 0$ in $\land V\otimes \land U$. Now a straightforward calculation shows that the connecting homomorphism is given explicitly by 
 \begin{eqnarray}
 \label{5}
 H(\varepsilon\otimes id): H^{\geq 1} (\land W\otimes \land U) \stackrel{\cong}{\longrightarrow} \mathbb Q a \otimes \land U.
 \end{eqnarray}

 On the other hand, applying  $-\otimes_{\land V}\land V\otimes \land U$ to the quasi-isomorphism $\gamma$ yields quasi-isomorphisms $
 \xymatrix{
 (\land Z, \overline{d}) &&   \land V\otimes \land U\otimes \land Z\ar[ll]^\simeq \ar[rr]^\simeq_\gamma && \land W\otimes \land U}$, so that we  have a degree 1 isomorphism $H^{\geq 1}(\land Z, \overline{d})   \stackrel{\cong}{\longrightarrow} \mathbb Q a \otimes \land U$. 
 It is immediate that $H(\varepsilon \otimes id)$ vanishes on products, which gives the second assertion.

 \hfill$\square$

 \vspace{3mm}
  Theorem 1 is now contained in

 \vspace{3mm}\noindent {\bf Theorem 1'.} {\sl Suppose $X$ is a connected CW complex, and $[f]\in \pi_n(X)$, some $n\geq 1$. Then in the factorization (3)  
 $$\xymatrix{\lambda : (\land V,d) \ar[r]^\eta & (\land V\otimes \land Z,d) \ar[r]^-\gamma_-\simeq & (\land W,d)},$$
the following conditions are equivalent:
\begin{enumerate}
\item[(i)] $[f]$ is rationally inert.
\item[(ii)] The generating space $Z$     can be chosen so that 
 \begin{eqnarray}
 \label{3}
 \overline{d}: Z\to \land^2Z  \hspace{5mm}\mbox{and }  H(\land Z) = \mathbb Q \oplus (Z\cap \mbox{ker}\, \overline{d}).
 \end{eqnarray}  
 \item[(iii)] The homotopy fibre of $F(f)$ of $i_{\mathbb Q}: X_{\mathbb Q}\to (X\cup_fD^{n+1})_{\mathbb Q}$ is rationally wedge-like.
 \end{enumerate}}

 \vspace{3mm}\noindent {\sl proof:}  \emph{(i) $\Longrightarrow$ (ii):} Since $\langle \lambda \rangle$ is identified with $i_{\mathbb Q}$,   $[f]\in\pi_n(X)$ is rationally inert if and only if the generating space $W$ can be chosen so that $\lambda$ restricts to an inclusion $V\to W$.
 In this case,  $\land W $ decomposes as a Sullivan  extension $\land V \to \land V\otimes \land Z = \land W$. Thus we may take $\eta= \lambda$ and $\gamma=id_{\land W}$. Note that if $\land V\otimes \land U$ is the acyclic closure of $\land V$, then the augmentation $\land V\otimes \land U\stackrel{\simeq}{\to}\mathbb Q$ defines a quasi-isomorphism $\land W\otimes \land U = \land V\otimes \land Z\otimes \land U \stackrel{\simeq}{\to} \land Z$.

 If dim$\, H^{\geq 1}(\land Z)= 1$, then necessarily $\land Z$ is the minimal Sullivan model of a sphere $S^k$ and $\langle \land Z\rangle = S^k_{\mathbb Q}$.  
 If dim$\, H^{\geq 1}(\land Z)\geq 2$,    let $\sigma : \land Z\to \land W\otimes \land U$ be a right inverse to the quasi-isomorphism $\land W\otimes \land U \simeq \land Z$ above. Since   $\land V\otimes \land Z$ is a  minimal Sullivan algebra,   it will follow that
 \begin{eqnarray}
 \label{7}
 \sigma : Z\to \land^{\geq 1}W\otimes \land U.
 \end{eqnarray}
 But this will imply that $\sigma: \land^{\geq 2}Z\to \land^{\geq 2}W\otimes \land U$. Now a simple calculation   shows that the connecting homomorphism vanishes on any $\overline{d}$-cycle in $\land^{\geq 2}Z$. Since the connecting homomorphism is an isomorphism it follows that division by $\land^{\geq 2}Z$ induces an injection $H^{\geq 1}(\land Z) \to Z$, and (ii) follows from Lemma 1.
 
 To complete this direction of the proof we need to establish (\ref{7}). For this write $Z = \cup_k Z(k)$ in which $Z(0)= Z\cap \mbox{ker}\, \overline{d}$ and $Z(k+1)= Z\cap \overline{d}^{-1}(\land Z(k))$. Assuming by induction that $\sigma : Z(k)\to \land^{\geq 1}W\otimes \land U$ we obtain that for $z\in Z(k+1)$, $d\sigma (z)=\sigma (\overline{d}z) \in \land^{\geq 2}W\otimes \land U$. Now let $\Phi$ be the component of $\sigma (z)$ in $1\otimes \land U$. Since $\land W$ is minimal it follows that $d : \land^{\geq 1}W\otimes \land U\to \land^{\geq 2}W\otimes \land U$. But if $\Phi \neq 0$ then $d(1\otimes \Phi)$ has a non-zero component in $V\otimes \land U$. Therefore $\Phi = 0$ and (\ref{7}) follows by induction on $k$.
 
 \vspace{3mm} \emph{(ii) $\Rightarrow$ (iii):}  
 Since $\langle \land Z\rangle \simeq F(f)$ it follows from Proposition 1 that $F(f)$ is rationally wedge-like.

 \vspace{3mm}\emph{(iii) $\Rightarrow$ (i):} First suppose that $F(f)$ is a rational sphere $S^k_{\mathbb Q}$. Then $\land Z$ is the minimal Sullivan model of a sphere, and so dim$\, Z\cap \mbox{ker}\, \overline{d}= 1$. Thus it follows from Proposition 4 that $U=0=V$. Since $\land V$ is the minimal Sullivan model for $X\cup_fD^{n+1}$ this implies that $\pi_*(X\cup_fD^{n+1})_{\mathbb Q} = 0$ and $[f]$ is rationally inert.
 
 Otherwise $F(f)$ is the inverse limit of rational wedges of at least two spheres. If $[f]$ is not inert then in the sequence
 $$\pi_*(\Omega (X\cup_fD^{n+1})_{\mathbb Q}) \to \pi_{*+1} (F(f)) \to \pi_{*+1} (X_{\mathbb Q})$$
 the image of $\pi_*(\Omega (X\cup_fD^{n+1})_{\mathbb Q}$ contains a non-zero class $\omega \in \pi_{*+1}(F(f))$. Because $\Omega (X\cup_fD^{n+1})_{\mathbb Q}$ acts on $F(f)$, it follows that the Whitehead product $\omega \bullet \beta$ of $\omega$ and any $\beta\in \pi_*(F(f))$ is zero.
 
 Then, because $\pi_*(F(f)) = \varprojlim \pi_*(S^{\sigma_1}\vee \dots \vee S^{\sigma_k})$ it follows that for some $r\geq 2$, the image $\overline{\omega}$ of $\omega$ in some $\pi_*(S^{\sigma_1}\vee\dots \vee S^{\sigma_k})_{\mathbb Q}$ is non-zero, and that
 $$\overline{\omega}\bullet \beta = 0,  \hspace{5mm} \beta\in \pi_*(S^{\sigma_1}\vee \dots \vee S^{\sigma_k})_{\mathbb Q}.$$
As observed in  (\ref{alpha}),  $\pi_*(S^{\sigma_1}\vee \dots \vee S^{\sigma_r})_{\mathbb Q}$ is the suspension of its homotopy Lie algebra $L$, and it follows from \cite[Chapter 2]{RHTII} that $\overline{\omega}$ determines a non-zero element in the center of $L$. But    the center of $L$ is zero, and therefore   $[f]$ is rationally inert.
 
 \hfill$\square$

 \section{Poincar\'e duality complexes}  We say a CW complex $Y= X\cup_fD^{n+1}$ is a \emph{rational Poincar\'e duality complex} if $H(Y)$ is a Poincar\'e duality algebra and the top class is in the image of $H(Y,X)$. In this case it follows that $H^{\leq n}(X) \stackrel{\cong}{\to} H(X)$. Poincar\'e duality complexes are rational Poincar\'e duality complexes, and so Theorem 2 follows from 
 
 \vspace{3mm}\noindent {\bf Theorem 2'.} {\sl If $Y = X\cup_fD^{n+1}$ is a rational Poincar\'e duality complex and the algebra $H(Y)$ requires at least two generators,  then $[f]\in \pi_n(X)$ is rationally inert.}
 
 \vspace{3mm} Before undertaking the proof we establish some notation. Let $\land V$ be the minimal Sullivan model of 
 $Y$, and let $S$ be a direct summand in $(\land V)^{n+1}$ of $(\land V)^{n+1}\cap \mbox{ker}\, d$. Then division by $S$ and by $(\land V)^{>n+1}$ defines a surjective quasi-isomorphism $\land V\stackrel{\simeq}{\to} A$, and 
 $$A^{n+1}= A^{n+1} \cap \mbox{Im}\, d \oplus \mathbb Q \omega,$$
 where $\omega$ is a cycle representing the top cohomology class of $Y$. As shown in (\cite[\S 5]{HL}), a cdga model of the inclusion $X\hookrightarrow Y$ is then provided by the inclusion
 $$j : (A,d) \to (A\oplus \mathbb Q t,d),$$
 where deg$\, t= n$, $t\cdot A^+=0$, and $dt= \omega$.
 
 Thus if $A\otimes \land U$ is the acyclic closure of $A$, then a cdga model for the homotopy fibre of $j$ is given by
 $$(A\oplus \mathbb Qt)\otimes_A(A\otimes \land U)= (A\oplus \mathbb Qt)\otimes \land U.$$
 Thus from the short exact sequence
 $$0\to A\otimes \land U \to (A\oplus \mathbb Qt)\otimes \land U\to \mathbb Q t\otimes \land U\to 0$$
 we deduce that
 $$H^{\geq 1}((A\oplus \mathbb Qt)\otimes \land U ) \stackrel{\cong}{\longrightarrow} \mathbb Qt\otimes \land U $$
 is an isomorphism of graded vector spaces.
 
 For the proof of Theorem 2 we first eliminate two special cases. First if $V^1= 0$ the argument of (\cite[\S 5]{HL}) shows that $(A\oplus \mathbb Qt)\otimes \land U$ is a cdga model of a wedge of spheres, and so $[f]$ is rationally inert. (Note that in \cite{HL} it is assumed that $X$ is simply connected; however the proof of this assertion relies only on the fact that $V^1= 0$.) Secondly, if $n=1$ then $X \simeq_{\mathbb Q} S^1_1\vee \dots \vee S^1_{2q}$ and so $Y$ is rationally equivalent to an oriented Riemann surface. In this case Theorem 2' is established in \cite{Artin}.

 Thus to prove Theorem 2' we may assume that $n\geq 2$ and that $A^1$ contains a non-zero cycle $x$. Since $H(A)$ is a Poincar\'e duality algebra there is a cycle $w\in A^n$ such that $wx=\omega$. The first step for the proof is then
 
 \vspace{3mm}\noindent {\bf Lemma 2.} {\sl With the hypotheses and notation above, $A^{n+1}\otimes \land U \subset d(A^n\otimes \land U)$.}
 
 \vspace{3mm}\noindent {\sl proof:} Choose $\overline{x}\in U^0$ so that $d\overline{x}=x$. Since $\land V$ is a minimal Sullivan algebra, $V$ is the union of an increasing sequence of subspaces $V(0)\subset \dots \subset V(q)\subset \dots$ in which $V(0)= \mathbb Q x$ and $d: V(q+1)\to \land V(q)$. It follows that $U$ is the union of an increasing sequence of subspaces $U(0)\subset \dots \subset U(q)\subset \dots$ in which $U(0)= \mathbb Q \overline{x}$ and $$d: U(q+1)\to A^{\geq 1}\otimes \land U(q).$$
 We show by induction on $q$ that
 \begin{eqnarray}
 \label{x}
 A^{n+1}\otimes \land U(q) \subset d(A^n\otimes \land U(q))
 \end{eqnarray}
 
 First note that any $z\in A^{n+1}$ has the form $z= dy + \lambda wx$, some $\lambda \in \mathbb Q$. Thus
 $$z\otimes 1 = d(y\otimes 1) \pm d(\lambda w\overline{x}) \in d(A^m\otimes \land U(0)).$$
 Then for $r\geq 1$,
 $$z\otimes \overline{x}^r = d(y\otimes \overline{x}^r \pm \frac{1}{r+1} w \otimes \overline{x}^{r+1}) + ry\otimes \overline{x}^{r-1}.$$
 It follows by induction on $r$ that $A^{n+1}\otimes \land U(0)\subset d(A^n\otimes \land U(0)).$
 
 Now fix a direct summand, $T$, of $U(q)$ in $U(q+1)$, and assume by induction that for some $s$,
 $$A^{n+1} \otimes \land U(q)\otimes \land^{\leq s}T \subset d(A^n\otimes \land U(q)\otimes \land ^{\leq s} T).$$ Then write $\Phi\in A^{n+1}\otimes \land U(q)\otimes \land^{\leq s+1}T$ as
 $\Phi = \sum \Phi_i \otimes \Psi_i$ with $\Phi_i \in A^{n+1}\otimes \land U(q)$ and $\Psi_i\in \land^{\leq s+1}T$. By the hypothesis $\Phi_i = d\Omega_i$ with $\Omega_i \in A^n\otimes \land U(q)$. Therefore
 $$\sum \Phi_i\otimes \Psi_i = d(\sum \Omega_i\otimes \Psi_i) \pm \sum \Omega_i \land d\Psi_i.$$
 The first term is in $d(A^n\otimes \land U(q)\otimes \land^{\leq s+1}T)$. On the other hand, $d\Psi_i \in A^{\geq 1} \otimes \land U(q)\land^{\leq s}T$ and so the second term is in $ A^{n+1}\otimes \land U(q)\otimes \land^{\leq s} T$. By hypothesis, the second term is contained in $d(A^n\otimes \land U(q)\otimes \land^{\leq s}T)$. This closes the induction. \hfill$\square$
 
 \vspace{3mm}\noindent {\sl proof of Theorem 2':} Let $\Phi\in \land U$. Then 
 $$t-(-1)^n w\overline{x} \in (A\oplus \mathbb Qt)\otimes \land U$$
 is a cycle, and 
 $$d((t-(-1)^nw\overline{x})\Phi) = -w\overline{x} \, d\Phi \in A^{n+1}\otimes \land U.$$
 By Lemma 2, $w\overline{x}\, d\Phi = d\Psi$ for some $\Psi \in A^n\otimes \land U$. Thus $(t-(-1)^n)w\overline{x}) \Phi + \Psi$ is a cycle projecting to $t\otimes \Phi$ in $\mathbb Q t\otimes \land U$. Then such cycles map to a basis of $\mathbb Q t\otimes \land U$. But because $n\geq 2$, $2n>n+1$ and so the product of any two of those cycles is zero. Therefore   this defines a cdga quasi-isomorphism from the cohomology of a wedge of spheres to $(A\oplus \mathbb Qt)\otimes \land U$. Lemma 1 and Theorem 1' together then imply that $[f]$ is rationally inert. \hfill$\square$

 \section{The structure of $L_Z$ and Theorem 3}
 
Any minimal Sullivan algebra $\land V$  equips $L_V$ with a natural additional structure (\cite[\S 3]{completions}), defined as follows. Associated with $\land V$ is the set, directed by inclusion, of the finite dimensional subspaces $V_\alpha\subset V$ for which $\land V_\alpha$ is preserved by $d$. For convenience we denote this set by
 $ {\mathcal J}_V = \{\alpha \}$. In particular,   
$$L_V = \varprojlim_{\alpha\in {\mathcal I}_V} L_\alpha, \hspace{5mm} L_\alpha \,\,\, \mbox{the homotopy Lie algebra of } \land V_\alpha.$$
That structure permits the explicit description of the Whitehead products in $\pi_*\langle \land V\rangle$ in terms of the Lie brackets in $L_V$ (\cite[Formula (11)]{completions}).  
 
 Moreover, for any augmented graded algebra, $A$, the \emph{classical completion} is defined by $\widehat{A}= \varprojlim_n A/I^n$, $I^n$ denoting the $n^{th}$ power of the augmentation ideal. The Sullivan completion of $UL_V$ is then the inverse limit, $$\overline{UL_V} = \varprojlim_\alpha \widehat{UL_\alpha}.$$
  Further,  by (\cite[Proposition 3.3]{depthLS1}), there are natural isomorphisms $\widehat{H}_*(\Omega \langle \land V_\alpha\rangle;\mathbb Q) \stackrel{\cong}{\longrightarrow} \widehat{UL_\alpha}$. Passing to inverse limits then yields the isomorphism of the \emph{Sullivan completions},
\begin{eqnarray}
\label{nine}
\overline{H_*(\Omega\langle \land V\rangle;\mathbb Q)} \stackrel{\cong}{\longrightarrow} \overline{UL_V}.
\end{eqnarray}

Similarly,  the \emph{Sullivan central series} is the filtration of $L_V$ given by
$$L_V^{(r)}= \varprojlim_{\alpha\in {\mathcal J}_V} L_\alpha^r,$$
where $L_\alpha^r$ is the ideal spanned by iterated commutators of length $r$. It satisfies  (\cite[\S 6]{completions})
$$L_V/L_V^{(r)} = \varprojlim_{\alpha \in {\mathcal J}_V} L_\alpha /L_\alpha^{(r)} \hspace{5mm}\mbox{and } L_V = \varprojlim_r L_V/L_V^{(r)}.$$
 In the case that $\langle \land V\rangle $ is the homotopy fibre of $i_{\mathbb Q} : X_{\mathbb Q}\to (X\cup_fD^{n+1})_{\mathbb Q}$ when $[f]\in \pi_n(X)$ is rationally inert, this additional structure has the striking properties provided in Theorem 3' below.

 Suppose next that $\land W = \land V\otimes \land Z$ is the decomposition of a minimal Sullivan  algebra determined   by an inclusion $\land V\to \land W$ with $V\subset W$, and denote $\mathbb Q\otimes_{\land V} \land W = (\land Z, \overline{d})$. Then the short exact sequence $V\to W\to Z$ dualizes to the short exact sequence
 $$0\leftarrow L_V\leftarrow L_W\leftarrow L_Z\leftarrow 0$$
 of Lie algebra morphisms, which   identifies $L_Z$ as an ideal in $L_W$. The \emph{holonomy representation} $\overline{\theta}$ of $L_V$ in $H(\land Z)$,   (\cite[Chapter 4]{RHTII}), then extends (\cite[\S 7]{completions}) to a     \emph{holonomy representation} of $\overline{UL_V}$ in $H(\land Z)$.

On the other hand,    the right adjoint representation of $L_{W}$ in $L_{Z}$   extends to  the right \emph{adjoint representation} of $\overline{UL_W}$ in $L_Z$, which further factors to give a right   representation   of $\overline{UL_V}$ in $L_Z/L_Z^{(2)}$ (\cite[Proposition 7]{completions}).

\vspace{3mm} Now suppose $(\land Z, \overline{d}$ is a quadratic Sullivan algebra. The surjection $\land^{\geq 1}Z\to Z$ with kernel $\land^{\geq 2}Z$ induces a surjection $H^{\geq 1}(\land Z) \to Z\cap \mbox{ker}\, \overline{d}$ of $\overline{UL_V}$-modules. This in turn dualizes to an inclusion
$$(Z\cap \mbox{ker}\, \overline{d})^\vee \to H^{\geq 1}(\land Z)^\vee$$
of right $\overline{UL_V}$-modules. Moreover, according to (\cite[Propositions 6 and 7]{completions}) the pairing $Z\times sL_Z \to \mathbb Q$ induces an isomorphism
\begin{eqnarray}
\label{douze}
L_Z/L_Z^{(2)} \stackrel{\cong}{\longrightarrow} (Z\cap \mbox{ker}\, \overline{d})^\vee
\end{eqnarray}
of right $\overline{UL_V}$-modules.

\vspace{2mm}
\emph{For the rest of this section we fix a map to a connected CW complex,
$$f: S^n \to X,$$
some $n\geq 1$, for which $[f]$ is rationally inert.}

\vspace{2mm} As observed in the Remark in $\S 1$, a Sullivan representative $\land V\to \land W$ for the inclusion $X\to X\cup_fD^{n+1}$ has the form 
$$\land V\to \land V\otimes \land Z= \land W,$$
and as above we denote the quotient differential in $\land Z$ by $(\land Z, \overline{d})$. It follows from Theorem 1' that $(\land Z, \overline{d})$ is a quadratic Sullivan algebra and that $H^{\geq 1}(\land Z, \overline{d})= Z\cap \mbox{ker}\, \overline{d}$. 

Now recall from $\S 2$ the linear map  
$$\varepsilon : \land W\to \mathbb Q$$
of degree $-n$. Since $\varepsilon (V)= 0$, $\varepsilon$ factors to give
$$\widehat{\varepsilon} \in (Z^n)^\vee = (L_Z)_{n-1}.$$
Thus, in view of (\ref{nine}), Theorem 3 is contained in 

\vspace{3mm}\noindent {\bf Theorem 3'.} {\sl With the hypotheses and notation above, let $\overline{\varepsilon}\in L_Z/L_Z^{(2)}$ denote the image of $\widehat{\varepsilon}$. Then 
\begin{enumerate}
\item[(i)] Both $L_Z/L_Z^{(2)}$ and $H^{\geq 1}(\land Z)^\vee$ are free $\overline{UL_V}$-modules, respectively generated by $\overline{\varepsilon}$ and $\widehat{\varepsilon}$. 
\item[(ii)] The map $\Phi \mapsto \varepsilon\cdot \Phi, \Phi \in \overline{UL_W}$, is a surjection
$$\tau : \overline{UL_W} \twoheadrightarrow L_Z,$$
of $\overline{UL_W}$-modules.
\item[(iii)] Any  subspace $S\subset L_Z$ with $S\stackrel{\cong}{\to} L_Z/L_Z^{(2)}$ freely generates a free   sub Lie algebra, $E\subset L_Z$, and
$$\varprojlim E/E\cap L_Z^{(r)} \stackrel{\cong}{\longrightarrow} L_Z.$$
\end{enumerate}}

\vspace{3mm}\noindent {\bf Remark.} When $X$ is simply connected with finite Betti numbers and $n\geq 2$, then Theorem 3'  is established in (\cite[Theorem 3.3]{HL}).

\vspace{3mm} Before undertaking the proof of Theorem 3' we establish a preliminary Proposition. For this, denote by $\varepsilon_W : \land V\otimes \land U \stackrel{\simeq}{\to} \mathbb Q$ the augmentation in the acyclic closure of $\land V$  defined by $\varepsilon_W(U)= 0$. Since the quotient differential in $\land U$ is zero, the holonomy representation of $\overline{UL_V}$ is a representation in $\land U$. On the other hand, the holonomy representation of $\overline{UL_V}$ in $H^{\geq 1}(\land Z)$ is a representation in $Z\cap \mbox{Ker}\, \overline{d}$. Now we strengthen Proposition 4 with

 \vspace{3mm}\noindent {\bf Proposition 5.} With the hypotheses and notation above, there is a commutative diagram
 $$\xymatrix{\land U \ar[rd]_{\varepsilon_W} \ar[rr]_\cong^\psi && Z\cap \mbox{Ker}\, \overline{d} 
 \ar[ld]^{\widehat{\varepsilon}}
 \\
 & \mathbb Q}$$
 in which $\psi$ is an isomorphism of $\overline{UL_V}$-modules of degree $n+1$.
 
 \vspace{3mm}\noindent {\sl proof.} Implicit in the isomorphism $\land W=\land V\otimes \land Z$ is the choice of a left inverse, $\land Z\to \land W$, of graded algebras  for the surjection $\land W\to \land Z= \mathbb Q\otimes_{\land V}\land W$. This, with $id_{\land V}$, defines an isomorphism $\land V\otimes \land Z \stackrel{\cong}{\to} \land W$, and identifies $id\otimes \widehat{\varepsilon}$ with $\varepsilon$. 
 A simple and standard argument  using Proposition 1  shows that this left inverse can be chosen so that the image of $\land V \otimes \left( (Z\cap \mbox{Ker}\, \overline{d})\oplus \mathbb Q\right)$  is preserved by $d$. 
 It is then immediate that the inclusion of this subcomplex in $(\land V\otimes \land Z)$ is a quasi-isomorphism. Thus from the commutative diagram (3) we obtain the   row exact sequence
 $$\xymatrix{
 0\to \mathbb Q a \ar[r] & \land V \otimes (Z\cap \mbox{Ker}\, \overline{d}\oplus \mathbb Q ) \oplus \mathbb Q a \ar[r] &   \land V \otimes (Z\cap \mbox{Ker}\, \overline{d}\oplus \mathbb Q) \ar[r] & 0\\
 \mbox{} & \land V.\ar[u]^\simeq \ar[ru]_\lambda}$$
 Since $\varepsilon (\land V)= 0$, $\land V$ is a   subcomplex. Division by this subcomplex   yields the row exact sequence of complexes,
 $$0 \to \mathbb Q a \to \land V\otimes (Z\cap \mbox{Ker}\, \overline{d}) \oplus \mathbb Q a \to \land V \otimes (Z\cap \mbox{Ker}\, \overline{d})\to 0$$
 in which the middle complex has zero homology. It is immediate that the connecting quasi-isomorphism $\delta$, is then given by
 $$\Phi\otimes z \mapsto \left\{\begin{array}{ll} \widehat{\varepsilon} (z)\, a &\mbox{if }\Phi = 1\\0 & \mbox{if }\Phi\in \land^{\geq 1}V.\end{array}\right.$$
 
 With a shift of degrees, regard $\varepsilon_W$ as a quasi-isomorphism $\land V\otimes \land U\stackrel{\simeq}{\to} \mathbb Q a$, sending $1\mapsto a$. Then, since $\land V\otimes \land U$ is $\land V$-semifree,  in the diagram,  
 $$\xymatrix{
 \land V\otimes \land U \ar@{-->}[rr]^\chi \ar[rd]^\simeq_{\varepsilon_W} && \land V\otimes (Z\cap \mbox{Ker}\, \overline{d})\ar[ld]_\simeq ^\delta \\
 & \mathbb Q a,}$$
 we may lift $\varepsilon_W$ through $\delta$ to obtain the quasi-isomorphism, $\chi$, of $\land V$-modules. But $\land V\otimes (Z\cap \mbox{Ker}\, \overline{d})$ is also $\land V$-semifree. Therefore applying $\mathbb Q \otimes_{\land V}-$ yields a quasi-isomorphism $\psi : \land U\stackrel{\simeq}{\to} Z\cap \mbox{Ker}\, \overline{d}$. 
 
Now the differentials in $\land U$ and in $Z\cap \mbox{Ker}\, \overline{d}$ are zero, and so $\psi$ is an isomorphism. Moreover,  $\mathbb Q \otimes_{\land V}-$ converts morphisms between $\land V$-semifree modules to morphisms of $L_V$-modules. In this case $\psi$ is then automatically a morphism of $\overline{UL_V}$-modules. Finally, it is also immediate that the diagram of the Proposition commutes. \hfill$\square$
 
 \vspace{3mm}\noindent {\sl proof of Theorem 2 (i).}  Here we rely consistently on the notation and conventions of $\S 2$. 
 
 First, observe that the dual of a $\overline{UL_V}$-module inherits a right $\overline{UL_V}$-module structure in the standard way. Thus replacing $\psi$ by $\psi^{-1}$ in the diagram of Proposition 5 and then dualizing yields the commutative diagram
\begin{eqnarray}
\label{8}\xymatrix{ (\land U)^\vee \ar[rr]^\cong &&  (Z\cap \mbox{Ker}\, \overline{d})^\vee\\
& \mathbb Q \ar[lu]\ar[ru]& ,}
\end{eqnarray}
in which $1\in \mathbb Q$ maps to $\varepsilon_W \in (\land U)^\vee$ and to $\widehat{\varepsilon}\in (Z\cap \mbox{ker}\, \overline{d})^\vee$. By (\cite[Proposition 8]{completions}) $(\land U)^\vee$ is a free right $\overline{UL_V}$-module, freely generated by $\varepsilon_W$. Since $H^{\geq 1}(\land Z)= (Z\cap \mbox{ker}\, \overline{d})$, it follows from (\ref{8}) that $H^{\geq 1}(\land Z)^\vee$ is a free right $\overline{UL_V}$-module freely generated by $\widehat{\varepsilon}$.

\vspace{2mm} (ii) To establish that the map
$$\tau : \overline{UL_W}\to L_Z$$
  is surjective, note that if $\beta\geq \alpha\in {\mathcal J}$ and $s\geq r$, then since $Z_\beta \supset Z_\alpha$,
$$L_{Z_\beta}/L_{Z_\beta}^s \longrightarrow L_{Z_\alpha}/L_{Z_\alpha}^r$$ is a surjection of finite dimensional spaces. Thus it is sufficient to show that the composites
\begin{eqnarray}
\label{(9)}
\overline{UL_W}\to L_Z\to L_{Z_\alpha}/L_{Z_\alpha}^{r+1}
\end{eqnarray}
are all surjective.

When $r= 1$, this is immediate from part (i) of the Theorem. Moreover, it follows from the construction of $\tau$ that its image is an ideal in $L_Z$. This, together with the surjectivity of  (\ref{(9)})  when $r= 1$ implies via the obvious induction that  (\ref{(9)})  is surjective for all $r$. 
 
 \vspace{3mm} (iii)  To show that $E$ is free it is sufficient to show that any linearly independent elements $x_1, \dots , x_k\in S$ generate a free sub Lie algebra $F$. But by (ii) the restriction of $S$ to $Z\cap \mbox{ker}\, \overline{d}$ is an isomorphism $sS \stackrel{\cong}{\to}  (Z\cap \mbox{ker}\, \overline{d})^\vee$. It follows that there are $z_1, \dots , z_k\in Z\cap \mbox{ker}\, \overline{d}$ such that
 $$\langle z_i, sx_j\rangle = \delta_{ij}.$$
 
 Let $T$ be the linear span of the $z_i$, so that $\mathbb Q\oplus T\subset \mathbb Q \oplus (Z\cap \mbox{ker}\, \overline{d})$ is a sub cdga, with minimal Sullivan model $\land Z_T\subset \land Z$ satisfying $T = Z_T \cap \mbox{ker}\, \overline{d}$, and with homotopy Lie algebra $L_T$. The surjection $L_Z\to L_T$ maps the generating set $\{x_i\}$ of $F$ bijectively to a dual basis for $T$. As shown in the Example in \S 2, it follows that $F$ is free.

 Finally, let $S_\alpha$ be the image of $S$ in $L_{Z_\alpha}$. Since $L_Z^{(2)} \to L_{Z_\alpha}^{(2)}$ is surjective, it follows that $S_\alpha +   L_{Z_\alpha}^{2} = L_{Z_\alpha}$. Therefore, because $L_{Z_\alpha}$ is nilpotent, the induced maps $E\to L_{Z_\alpha}$ are surjective. Hence, these induce surjections $E/E\cap L_Z^{(r)} \to L_{Z_\alpha}/L_{Z_\alpha}^r$. 
 
Since each $L_{Z_\alpha}/L^r_{Z_\alpha}$ is finite dimensional, it follows that passing to inverse limits yields surjections
 $$E/E\cap L_Z^{(r)} \to L_Z/L_Z^{(r)}.$$
 It is immediate from this that $\varprojlim_r E/E\cap L_Z^{(r)} \stackrel{\cong}{\longrightarrow} L_Z.$
 
 \hfill$\square$

  \section{One-relator groups}
  
  Our objective here is the proof of
  
  \vspace{3mm}\noindent {\bf Theorem 4} If $X$ is a wedge of at least two circles then any non-zero $[f]\in \pi_1(X)$     is rationally inert  or, equivalently, $(X\cup_fD^2)_{\mathbb Q}$ is aspherical.
  
  \vspace{3mm}\noindent {\sl proof:} First observe that in fact
  \begin{eqnarray}
  \label{13}
  \mbox{$[f]$ is rationally inert $\Leftrightarrow$ $(X\cup_fD^2)_{\mathbb Q}$ is aspherical.}
  \end{eqnarray}
  In fact, the same argument as in the Example in \S 2 shows that the minimal Sullivan model of $X$ is cdga equivalent to $\mathbb Q \oplus H^1(X_{\mathbb Q})$. It follows that the homotopy Lie algebra, $L$, is concentrated in degree $0$ and since $\pi_*(X_{\mathbb Q})= sL$, $X_{\mathbb Q}$ is aspherical. Thus if $[f]$ is rationally inert then $(X\cup_fD^2)_{\mathbb Q}$ is aspherical. On the other hand, a Sullivan representative for the inclusion $i : X\to X\cup_fD^2$ is a morphism $\gamma : \land V\to \land W$ of minimal Sullivan algebras. Since $\pi_1(i)$ is injective, $H^1(i)$ is surjective  and it follows that $\gamma : V^1\to W^1$ is injective. But if $(X\cup_fD^2)_{\mathbb Q}$ is aspherical, then $V = V^1$, $\gamma $ is injective, and by definition $[f] $ is rationally inert.
  
  Next note that it is sufficient to prove the Theorem when $X$ is a finite wedge of circles. Simply write $X = Y\vee Y'$ in which $Y$ is a finite wedge of circles, $Y'$ is a wedge of circles, and $f: S^1\to Y$. Then, as just observed,    $Y'_{\mathbb Q}$ is aspherical.  It follows from Proposition 2 that if $(Y\cup_fD^2)_{\mathbb Q}$ is aspherical, then so is $(X\cup_fD^2)_{\mathbb Q} = \left[ (Y\cup_fD^2)\vee Y'\right]_{\mathbb Q}$. Thus by (\ref{13}), $[f]\in \pi_1(Y\cup_fD^2)_{\mathbb Q}$ is rationally inert if and only if $[f]\in \pi_1(X_{\mathbb Q})$ is rationally inert.

 \vspace{2mm} 
  \emph{In summary, we may and do assume henceforth that}
  $$X = S^1\vee \dots \vee S^1.$$
  
  On the other hand, we observe that
\begin{eqnarray} [f]\neq 0 \hspace{3mm} \Rightarrow \hspace{3mm} \mbox{a Sullivan representative of $f$ is non-zero.}\end{eqnarray}
  In fact, denote $G = \pi_1(X)$, so that $G_{\mathbb Q}= \pi_1(X_{\mathbb Q})$. According to \cite[Theorem 7.5]{RHTII}, $G^n/G^{n+1} \otimes \mathbb Q \stackrel{\cong}{\to} G^n_{\mathbb Q}/G^{n+1}_{\mathbb Q}$. But by \cite{hall}, $G^n/G^{n+1}$ is a free abelian group, and hence $G^n/G^{n+1}\to G^n_{\mathbb Q}/G^{n+1}_{\mathbb Q}$ is injective. Since $G$ is a free group, $G\to \varprojlim_n (G/G^n)_{\mathbb Q}$ is injective and the image of $[f] $ in $G_{\mathbb Q}$ is non-zero.
  In particular, a Sullivan representative of $f$ is non-zero.
  
  Next recall from the Example in \S 2 and Lemma 1 that $S^1\vee \dots \vee S^1\vee S^2$ has a quadratic minimal Sullivan model, $(\land W, d_1)$ in which $W\cap \mbox{ker}\, d_1= H^{\geq 1} (S^1\vee \dots \vee S^1\vee S^2)$. In particular, $W^1\cap \mbox{ker}\, d_1= H^1(S^1\vee \dots \vee S^1)$.
  Moreover, $W^{>1}\cap \mbox{ker}\, d_1= W^2\cap \mbox{ker}\, d_1= \mathbb Q a$, where $a$ represents the orientation class of $S^2$. It follows that
  $$W = W^1\oplus \mathbb Q a \oplus R,$$
  and that the identity in $\land W^1$ extends to a quasi-isomorphism
  $$\varphi : (\land W,d_1) \stackrel{\simeq}{\to} (\land W^1\oplus \mathbb Q a, d_1)$$
  with $\varphi (a) = a$ and $\varphi (R)= 0$. 
  
  \vspace{3mm}\noindent {\bf Note:}  In comparing with the general situation described in \S 3, observe that the $\land W^1$ here corresponds to the $\land W$ in \S 3, and that the $\land W$ here has no analogue in \S 3.
  
  \vspace{3mm}
  In particular $\varphi$ preserves wedge degrees when $a$ is assigned wedge degree 1. Thus not only is $H(\mbox{ker}\, \varphi) = 0$, but in fact for cycles $\Phi\in \land W$, 

\begin{eqnarray}
\label{14}\Phi \in \land^kW\cap \mbox{ker}\, \varphi \Longrightarrow \Phi = d_1\Psi \hspace{3mm}\mbox{for some }\Psi\in \land^{k-1}W\cap \mbox{ker}\, \varphi\,.
  \end{eqnarray}

  The proof of Theorem 3 is now accomplished in the following steps:
  
  \vspace{2mm} \noindent \emph{Step One: Construction of a linear map of degree 1, $d_0 : W\to W$, whose extension, also denoted $d_0$, to a derivation in $\land W$ provides a cdga $(\land W, d_1+d_0)$ connected by cdga quasi-isomorphisms to $A_{PL}(X\cup_fD^2)$. } 
  
  \vspace{3mm} \noindent \emph{Step Two:  $(\land W, d_0+d_1)$ is a Sullivan algebra, and hence a Sullivan model for $X\cup_fD^2$.}
  
  \vspace{3mm}\noindent \emph{Step Three:  The minimal Sullivan model of $(\land W, d_1+d_0)$ has the form $(\land V^1,D)$, and so $(X\cup_fD^2)_{\mathbb Q}$ is aspherical, and $[f]$ is rationally inert.}

\vspace{3mm}\noindent \emph{Step One:    Construction of     $d_0: W\to W$ whose extension to a derivation (also denoted by $d_0$) provides a   cdga $(\land W, d_0+d_1)$ connected by cdga quasi-isomorphisms to  $A_{PL}(X\cup_fD^2)$.}

\vspace{1mm}  For this, fix   a   Sullivan representative $\psi : (\land W^1,d) \to (\land v,0)$ for $f$  and, as at the start of \S 3, define   $\varepsilon : \land W^1\to \mathbb Q$ by
$$\varepsilon (1)= \varepsilon (\land^{\geq 2}W^1)=0 \hspace{3mm}\mbox{and } \psi (w) = \varepsilon (w)v, \hspace{3mm} w\in W^1.$$
Then  define a derivation $\delta$ in $\land W^1\oplus \mathbb Q a$ by setting
  $$\delta (w) =\varepsilon (w) a\hspace{5mm}\mbox{and } \delta (\land^{\geq 2} W^1 \oplus \mathbb Q a)= 0.$$
  Then $d_1\delta = 0= \delta d_1$ and $\delta^2= 0$, so that $(\land W^1\oplus \mathbb Q a, d_1+\delta)$ is a cdga. As observed at the start of \S 3,   this cdga is connected by cdga quasi-isomorphisms to $A_{PL}(X\cup_fD^2)$.
  
Next, we construct a linear map $d_0 : W\to W$ of degree $1$ such that $d_0d_1+d_1d_0= 0$ and $\varphi \circ d_0= \delta \circ \varphi$.   
  
  For this, recall that $W=\cup_n W(n)$ with $ W(0) = W\cap \mbox{ker}\, d_1$ and $W(n+1) = W\cap d_1^{-1}(\land W(n))$.  By convention, $W(-1)= 0$.   We assume by induction that $d_0$ is constructed in $W(n-1)$, and write $W(n)= W(n-1)\oplus S$. If $w\in S$, then
  $$d_1d_0d_1w = -d_0d_1^2w = 0\,,$$
  and so $d_0d_1w$ is a cycle in $(\land^2W,d_1)$.
  
  Suppose first that $w\in W^1$. Then $d_1w\in \land^2W^1(n-1)$ and
  $$\varphi (d_0d_1w) = \delta \varphi (d_1w)= 0\,.$$
  Thus by (\ref{14}), for some $u\in  \mbox{ker}\, \varphi \cap W^2$,
  $$d_0d_1w = d_1u\,.$$
  Moreover, $\delta : W^1\to \mathbb Q a$, and so we may regard $\delta w$ as an element of $W^2$ for which $d_1\delta w= 0$ in $\land W$. Set $d_0w = \delta w-u$. Then
  $$d_1d_0w = -d_1u = -d_0d_1w$$
  and,
  since $\varphi u = 0$, $$\varphi (d_0w) = \varphi (\delta w) = \delta w= \delta (\varphi w)\,.$$
  
  On the other hand suppose $w\in W^k$, some $k\geq 2$. Then $d_0d_1w\in (\land^2W)^{k+2}$ and so $d_0d_1w\in R\land \land W \oplus \mathbb Q a^2$. Thus $\varphi (d_0d_1w)= 0$ and again by (\ref{14}) $d_0d_1w=d_1u$ for some $u\in W^{\geq 3}\subset R$. Set $d_0w= -u$, so that again
  $$d_1d_0w = -d_1u = -d_0d_1w.$$
 Then, since $u\in R$, $\varphi u= 0$ while $\varphi w \in \mathbb Q a$ and so $\delta \varphi w= 0$ as well. This completes the construction of $d_0$. By construction,
  $$\varphi \circ (d_1+d_0) = (d_1+\delta) \circ \varphi\,.$$
  
Finally we show that $d_0^2 = 0$ so that $d_1+d_0$ is a differential, and that
  \begin{eqnarray}
  \label{(16)}
  \varphi : (\land W, d_1+d_0) \stackrel{\simeq}{\longrightarrow} (\land W^1\oplus \mathbb Q a, d_1+\delta).
  \end{eqnarray}
  In fact $d_1d_0^2 = d_0^2d_1$. Assume by induction that $d_0^2= 0$ in $W(n-1)$. Then for $w\in S$, $d_0^2w$ is a $d_1$-cycle and $\varphi (d_0^2w) = \delta^2\varphi w = 0$. Thus by (\ref{14}), $d_0^2w$ is a $d_1$-boundary, and hence $d_0^2w= 0$. Thus $(\land W, d_1+d_0)$ is a cdga and  $\varphi$ is a morphism of cdga's with respect to $d_1+d_0$ and $d_1+\delta$. Filter both sides by the difference between degree and wedge degree. The map induced by $\varphi$ in the $0^{th}$ term of the spectral sequence is the quasi-isomorphism $\varphi : (\land W, d_1) \stackrel{\simeq}{\longrightarrow} (\land W^1\oplus \mathbb Q a, d_1)$. This establishes (\ref{(16)}) 
  
  Note that by (16), the Sullivan representative $\psi$ is non-zero, and so for some $w\in W^1$, $\delta w= a$, and $d_0w\neq 0$.
  
 \vspace{3mm} \emph{Step Two:  $(\land W, d_1+d_0)$ is a Sullivan algebra, and hence is a Sullivan model for $X\cup_fD^2$.}
  
  \vspace{2mm} Here we prove a more general result: if   $(\land V,d)$ is any minimal Sullivan algebra and $d_0: V\to V$ is a linear map of degree 1 such that $d_0^2= dd_0+d_0d= 0$, then $(\land V, d+d_0)$ is a Sullivan algebra.

  For this, fix an increasing filtration $0= V(0)\subset \dots \subset V(n)\subset \dots $ such that $V = \cup_n V(n)$ and $d : V(n+1)\to \land^{\geq 2}V(n)$. Then, as follows,  define by induction a sequence of subspaces of $V$ of the form
  $$Q(0)\subset P(0)\subset \dots \subset Q(n)\subset P(n)\subset \dots$$
  so that
  $$d \hspace{1mm}\mbox{and }d_0 : Q(n+1) \to \land P(n),$$
  $$d\hspace{1mm}\mbox{and }d_0 : P(n+1)\to \land Q(n+1),$$
  and
  $$P(n)\supset V(n).$$
  
  First, we set $Q(0)= P(0) = 0$. Then suppose $Q(k)$, and $P(k)$ are constructed for $k\leq n$. Write
  $$V(n+1) = V(n+1) \cap P(n) \oplus S(n+1),$$
  and set 
  $$Q(n+1) = P(n) + d_0(S(n+1)) \hspace{5mm}\mbox{and } P(n+1) = Q(n+1) + S(n+1).$$
  It is immediate that
  $$P(n+1) \supset P(n) + S(n+1) \supset V(n) + S(n+1) = V(n+1).$$
  
  Moreover, if $x\in S(n+1)$ then
  $$\renewcommand{\arraystretch}{1.4}\begin{array}{ll}
  dd_0x &= -d_0dx \in d_0d(S(n+1))\subset d_0(\land^{\geq 2}V(n))\\
  & \subset d_0(\land^{\geq 2}P(n)) \subset \land^{\geq 2}P(n)\,.
  \end{array}\renewcommand{\arraystretch}{1}$$
  In particular, $d: Q(n+1) \to \land^{\geq 2}P(n)$. Further $d_0^2(S(n+1)) = 0$ and so $d_0(Q(n+1)) = d_0(P(n)) \subset P(n)$. 
  
  On the other hand, if $x\in S(n+1)$ then $d_0x\in Q(n+1)$ by construction, while $dx\in \land^{\geq 2}V(n)\subset \land P(n)$. This closes the induction and exhibits $(\land V, d+d_0)$ as a Sullivan algebra.

   \vspace{3mm}\noindent \emph{Step Three: The minimal Sullivan model of $(\land W, d_1+d_0)$ has the form $(\land V^1,D)$, and so $(X\cup_fD^2)_{\mathbb Q}$ is aspherical.}
   
  \vspace{1mm} Recall from the Example in \S 2  that the homotopy Lie algebra of $(\land W, d_1)$ is the completion, $\widehat{\mathbb L}$ of the free Lie algebra $\mathbb L (x_1, \dots , x_r,y)$ generated by vectors $x_i$ dual to the orientation classes of the circles, and by $y$ dual to the orientation class of $S^2$. By construction, $W^{\geq 2} = \mathbb Q a \oplus R$, and we may choose $y$ so that 
  $$\langle a, sy\rangle = 1 \hspace{3mm}\mbox{and } \langle R, sy\rangle = 0.$$
Now dualize $d_0 : W\to W$   to $d: \widehat{\mathbb L}\to \widehat{\mathbb L}$. Since deg$\, d= -1$ it follows that $d: \widehat{\mathbb L}(x_1, \dots , x_r)\to 0$ and $dy\in \widehat{\mathbb L}(x_1, \dots , x_r)$. Moreover, because $d_0$ is a derivation satisfying $d_0d_1+d_1d_0= 0= d_0^2$, it follows that $d$ is a derivation in the Lie algebra $\widehat{\mathbb L}$ and that $d^2= 0$. 
  
  Moreover, if $(\land V,D)$ is the minimal Sullivan model of $(\land W,d_1+d_0)$ then $V \cong H(W, d_0)$. Therefore  $H(\widehat{\mathbb L},d) =  (H(W, d_0))^\vee$, and so it is sufficient to prove that 
  $$H_{\geq 1}(\widehat{\mathbb L},d) = 0.$$

  Recall also from Step One that a Sullivan representative for $f$ determines a linear map $\varepsilon : W^1\to \mathbb Q$. Thus $\varepsilon $   desuspends to $\alpha \in  {L}_{W^1}= \widehat{\mathbb L}(x_1, \dots , x_r)$. We show now that
  \begin{eqnarray}
  \label{dixneuf}
  dy=\alpha,\end{eqnarray}
  so that $dy\neq 0$.
  
  For this, recall from Step One that if $w\in W^1$ then $d_0w = \varepsilon (w) a- u$, where $u\in W^2\cap \mbox{ker}\, \varphi= Z$. It follows that 
  $$\langle w,sdy\rangle = -\langle d_0w, sy\rangle = -\langle \varepsilon (w)a-u, sy\rangle = \langle w, s\alpha\rangle,$$
which establishes (\ref{dixneuf}). 
  
  Denote by $\mathbb L_q(x_i)$ the linear span of the commutators of length $q$ in the $x_i$. Write
  $dy$ as a series
  $$dy= \sum_{q\geq n} \alpha_q$$
  where $\alpha_q\in \mathbb L_q(x_i)$ and $\alpha_n \neq 0$. Then form the differential graded Lie algebra $(\mathbb L(x_i,y), \partial)$ with $\partial (x_i)= 0$ and $\partial (y) = \alpha_n$. Since $\alpha_n$ belongs to $\mathbb L_n(x_i)$ we can modify  the degrees in $\mathbb L(x_i)$ by assigning deg$\, 2$ to the $x_i$,   without changing the homology with respect to $\partial$. Thus it follows from \cite[Theorem 3.12]{HL} that   $H_q(\mathbb L(x_i, y), \partial)= 0$ for $q>0$. 
  
  Now let $\omega= \sum_{q\geq p}\omega_q$ be a $d$-cycle in degree $r>0$ in $\widehat{\mathbb L}(x_i, y)$, with $\omega_q\in \mathbb L_q(x_i,y)$. Then $\omega_p $ is a $\partial$-cycle, and so a $\partial$-boundary. Choose $\beta_{p-n+1}\in \mathbb L_{p-n+1}(x_i, y)$ with $\partial (\beta_{p-n+1})= \omega_p$. Write $\omega(1) = \omega - d(\beta_{p-n+1})$, then $\omega (1)$ is a sum $\sum_{s\geq p+1}\omega(1)_s$. One again $\omega(1)_{p+1}$ is a $\partial$-cycle. This determines $\beta_{p-n+2}$. Continue in this way to obtain at the and an element
$$\beta=  \sum_{s\geq p-n+1} \beta_s$$ with $d\beta = \omega$.

 \hfill$\square$

  \vspace{3mm}\noindent {\bf Corollary.}  With the notation of Theorem 3,  set $V = W^1\cap \mbox{ker}\, d_0$. Then $(\land V,d_1)\to (\land W^1, d_1+d_0) $ is the minimal Sullivan model of $X\cup_fD^2$.   
  
  \vspace{3mm}\noindent {\sl proof:}  First note that any element in $\land^2W^1$ can be written as $\Phi = \sum_{i=1}^n w_i\land w_i'$ in which $w_1, \dots , w_n, w_1', \dots , w_n'$ are all linearly independent. Thus if $d_0\Phi= 0$ then each $d_0w_i=d_0w_i'= 0$. But $d_1: V\to \land^2W^1\cap \mbox{ker}\, d_0$, and so $\land V$ is preserved by $d_1$. It is immediate from Step Three that $V\stackrel{\cong}{\to} H(W,d_0)$, and it follows that $(\land V, d_1)$ is the minimal Sullivan model of $X\cup_fD^2$.  
  
  \hfill$\square$
 
 \section{Whitehead's problem and Theorem 5}

 \vspace{3mm}\noindent {\bf Theorem 5.}  If $X$ is a connected CW complex and $(X\cup \amalg_{k=1}^p D^2)_{\mathbb Q}$ is aspherical then $X_{\mathbb Q}$ is aspherical. 
 
 \vspace{3mm}\noindent {\sl proof:} The obvious induction reduces the statement to the case $p=1$. Then, since $\pi_*((X\cup_fD^2)_{\mathbb Q}) \cong V^\vee$ as sets where $\land V$ is the minimal Sullivan model of $X\cup_fD^2$, our hypothesis simply implies that $V = V^1$. Let $\varphi : (\land V,d)\to (\land W,d)$ be a Sullivan representative for the inclusion $i : X \to X\cup_fD^2$.  Since $H^1(X\cup_fD^2)\to H^1(X)$ is injective, it follows that $\varphi$ is injective and so $\land W$ decomposes as $\land V\otimes \land Z$, with $Z= Z^{\geq 1}$. In particular $\vert f\vert$ is rationally inert. Moreover, it follows from Proposition 1 that 
 $$H^{\geq 1}(\land Z, \overline{d}) \cong \mathbb Q b\otimes \land U,$$
 where deg$\, b= 1$ and $\land V\otimes \land U$ is the acyclic closure of $\land V$. Since $V = V^1$, $U= U^0$ and $H^{\geq 1}(\land Z, \overline{d}) = H^1(\land Z, \overline{d})$. This in turn implies $Z= Z^1$ and $X_{\mathbb Q}$ is aspherical.
 
 \hfill $\square$

 \vspace{5mm}\noindent Institut de Math\'ematique et de Physique, Universit\'e Catholique de Louvain, 2, Chemin du cyclotron, 1348 Louvain-La-Neuve, Belgium, yves.felix@uclouvain.be

 \vspace{1mm}\noindent Department of Mathematics, Mathematics  Building, University of Maryland, College Park, MD 20742, United States, shalper@umd.edu

\end{document}